\documentclass{article}
\usepackage{amsmath,amssymb,hyperref}

\begin{document}

\title{Basel problem: a physicist's solution}
\author{Zurab K. Silagadze}
\date{}

\maketitle

To paraphrase renowned American physicist Richard Feynman's quote 
``every theoretical physicist who is any good knows six or seven different 
theoretical representations for exactly the same physics'' \cite{1}, every
mathematician who is any good knows a dozen of solutions of the Basel problem,
which asks for an evaluation of the infinite series $\sum\limits_{n=1}^\infty 
\frac{1}{n^2}$. For example, Moreno in the arXiv version of \cite{2} gives 
more than 80  references related to various proofs of Euler's famous formula, 
and some new ones have appeared since then \cite{3,4,4A,5,6,7,7A}. Even Euler 
himself gave at least four proofs \cite{8,9,9A} that
\begin{equation}
\zeta(2)\equiv\sum\limits_{n=1}^\infty \frac{1}{n^2}=\frac{\pi^2}{6}.
\label{eq1}
\end{equation}
In \cite{10} W\"{a}stlund reformulated the Basel problem in terms of 
a physical system using the proportionality of the apparent brightness of 
a star to the inverse square of its distance. Inspired by this approach, 
here we give another physical interpretation which, in our opinion, is 
simpler, natural enough, and leads to a proof of (\ref{eq1}) which is very 
Eulerian in its spirit.

\section*{A physicist's solution of the Basel problem}
First of all, let's notice \cite{10} that, because of
$$\sum\limits_{n=1}^\infty \frac{1}{n^2}=\sum\limits_{n=1}^\infty \frac{1}
{(2n-1)^2}+\frac{1}{4}\sum\limits_{n=1}^\infty \frac{1}{n^2},$$
we have
\begin{equation}
\sum\limits_{n=1}^\infty \frac{1}{n^2}=\frac{4}{3}\sum\limits_{n=1}^\infty 
\frac{1}{(2n-1)^2}=\frac{1}{3}\sum\limits_{n=1}^\infty \frac{1}{(n-\frac{1}
{2})^2}.
\label{eq2}
\end{equation}
A physicist can interpret (\ref{eq2}) as representing a Coulomb force 
exerted on a unit charge located at a point $x=1/2$ by a semi-infinite 
sequence of opposite charges located at points $x_1=1$, $x_2=2$, \dots 
The corresponding electrostatic potential is
\begin{equation}
U(x)=-\sum\limits_{n=1}^\infty \frac{1}{n-x},
\label{eq3}
\end{equation}
and we can write
\begin{equation}
\sum\limits_{n=1}^\infty \frac{1}{n^2}=\left . \frac{1}{3}F(x)\right |
_{x=\frac{1}{2}}=\frac{1}{3}\left(-\left . \frac{dU(x)}{dx}\right)
\right |_{x=\frac{1}{2}}.
\label{eq4}
\end{equation}
Unfortunately (\ref{eq3}) diverges and hence we are returning here 
(temporarily) to the standards of mathematical rigor of Eulerian times. 
However physicists are used to infinities, and thus let's regularize the
potential (\ref{eq3}):
\begin{equation}
U(x)\to U_R(x)=U(x)-U(0)=-\sum\limits_{n=1}^\infty \left (\frac{1}{n-x}-
\frac{1}{n}\right ).
\label{eq5}
\end{equation}
Note that our regularization procedure does not affect at all the force and
hence we have
\begin{equation}
F(x)=\frac{d}{dx}\sum\limits_{n=1}^\infty \left (\frac{1}{n-x}-
\frac{1}{n}\right ).
\label{eq6}
\end{equation}
In (\ref{eq6}) we recognize immediately the presence of the digamma function 
$\psi(x)=\Gamma^\prime(x)/\Gamma(x)$ because
\begin{equation}
\psi(1-x)=-\gamma-\sum\limits_{n=1}^\infty \left (\frac{1}{n-x}-
\frac{1}{n}\right ),
\label{eq7}
\end{equation}
$\gamma$ being the Euler constant (this relation will be discussed below). 
Therefore
$$F(x)=\frac{d}{dx}\left[-\psi(1-x)-\gamma\right]=\psi_1(1-x),$$
where $\psi_1(x)=\psi^\prime(x)$ is the trigamma function. Fortunately
(\ref{eq4}) indicates that we need the trigamma function at $x=1/2$:
\begin{equation}
\sum\limits_{n=1}^\infty \frac{1}{n^2}=\frac{1}{3}\,\psi_1\left(\frac{1}{2}
\right ),
\label{eq8}
\end{equation}
and this is just the value at which the trigamma function can be simply 
calculated thanks to the Euler's reflection formula:
\begin{equation}
\psi_1(1-x)+\psi_1(x)=\frac{\pi^2}{\sin^2{\pi x}},
\label{eq9}
\end{equation}
which gives $\psi_1\left(\frac{1}{2}\right )=\frac{\pi^2}{2}$, and hence from 
(\ref{eq8}) we immediately get Euler's famous formula (\ref{eq1}).

\section*{Several remarks}
Above we presented a physics-motivated approach to the Basel problem.
Of course the connection with physics is tenuous at best. However the 
interpretation of inverse squares as representing Coulomb forces was a crucial
insight in defining a connection with polygamma functions and the reflection
formula. The resulting formalism is, in fact, quite elementary, in the sense 
that its basic pillars (\ref{eq7}) and (\ref{eq9}) can be obtained by 
elementary means. 

For example, (\ref{eq7}) follows from Newman's infinite product formula 
(used by Weierstrass as his definition of the gamma function)
\begin{equation}
\frac{1}{\Gamma{(1+z)}}=e^{\gamma z}\prod\limits_{n=1}^\infty\left(1+\frac{z}
{n}\right)e^{-z/n},
\label{eq28}
\end{equation}
which by itself is just another version of Euler's definition of the gamma
function as the limit \cite{22} 
\begin{equation}
\Gamma(z)=\lim_{n\to\infty}\frac{n^zn!}{z(z+1)\cdots (z+n)}.
\label{eq29}
\end{equation}
A simple consequence of (\ref{eq29}) is the following interesting identity
\cite{23} 
\begin{equation}
\prod\limits_{k=0}^\infty\frac{(k+\alpha_1)\cdots(k+\alpha_n)}
{(k+\beta_1)\cdots(k+\beta_n)}=\frac{\Gamma(\beta_1)\cdots\Gamma(\beta_n)}
{\Gamma(\alpha_1)\cdots\Gamma(\alpha_n)},
\label{eq30}
\end{equation}
where $n\ge 1$ and $\alpha_1,\ldots,\alpha_n$, $\beta_1,\ldots,\beta_n$
are nonzero complex numbers, none of which are negative integers, such that
$\alpha_1+\ldots+\alpha_n=\beta_1+\cdots+\beta_n$. In particular, when 
$\alpha_1=1+z$, $\alpha_2=1-z$, $\beta_1=\beta_2=1$, we get
$$\prod\limits_{n=1}^\infty\left(1-\frac{z^2}{n^2}\right)=\frac{1}
{\Gamma{(1+z)}\,\Gamma{(1-z)}},$$
which, in combination with $\Gamma{(1+z)}=z\,\Gamma(z)$ and Euler's celebrated
formula
\begin{equation}
\prod\limits_{n=1}^\infty\left(1-\frac{z^2}{n^2}\right)=
\frac{\sin{\pi z}}{\pi z},
\label{eq31}
\end{equation}
implies the validity of the reflection formula
\begin{equation}
\Gamma(z)\,\Gamma{(1-z)}=\frac{\pi}{\sin{\pi z}},
\label{eq32}
\end{equation}
from which other reflection formulas, like (\ref{eq9}), do follow.

It is tempting to consider the infinite product formula (\ref{eq31}) to be 
a real backbone of the presented approach, as in the Euler's original first 
proof. Although several elementary proofs of Euler's infinite product for 
the sine exist in the literature (see, for example, \cite{24,25,26,27}), they 
do not seem to be significantly simpler than the original proof by Euler. 
Therefore it may appear that the reflection formula is fairly nontrivial to 
derive and its proof is as hard a problem as the one we seek to solve.
However this is actually not the case. It is possible to avoid the use of the 
Euler's infinite product in the derivation of the reflection formula. Below 
we provide one such proof which is simple and elementary enough and doesn't 
rely on the infinite product formula for the sine function.

This proof of the reflection formula was inspired by Richard Dedekind's 1852 
proof \cite{28} of (\ref{eq32}) which seems to be not as well known as it 
deserves to be. It appears as an exercise in \cite{29} and was popularized in 
\cite{30}. We prove not (\ref{eq32}), but the reflection formula for the 
digamma function
\begin{equation}
\psi(x)-\psi(1-x)=-\pi\cot{\pi x},
\label{eq33}
\end{equation}
from which (\ref{eq9}) follows by differentiation.

During the proof, which seems to be much simpler than the Dedekind's original 
one, we freely interchange the order of integrals and differentiate under the 
integral signs, as physicist are generally accustomed to doing. A genuine 
mathematician, of course, will resort in these cases to Fubini's theorem 
and to Lebesgue's dominated convergence theorem to justify these 
operations \cite{30}.

We have the following well known integral representation for the digamma 
function:
\begin{equation}
\psi(1-x)=-\gamma+\int_0^1\frac{1-t^{-x}}{1-t}\,dt.
\label{eq34}
\end{equation}
Indeed, expanding $1/(1-t)$ in geometric series, interchanging the order of
summation and integration and thus integrating term by term, we get 
(\ref{eq7}).  

Therefore,
\begin{equation}
\phi(x)\equiv\psi(x)-\psi(1-x)=\int_0^1\frac{t^{-x}-t^{x-1}}{1-t}\,dt=
\frac{1}{2}\int_0^\infty\frac{t^{-x}-t^{x-1}}{1-t}\,dt.
\label{eq35}
\end{equation}
The validity of the last step can be checked by breaking the corresponding 
integral into two integrals, over $(0,1)$ and $(1,\infty)$, and putting
$y=1/t$ in the second integral.

Naively, the two parts of the last integral in (\ref{eq35}) 
\begin{equation}
\int_0^\infty\frac{t^{-x}}{1-t}\,dt\;\;\mathrm{and}\;\;
\int_0^\infty\frac{-t^{x-1}}{1-t}\,dt,
\label{eq36}
\end{equation}
appear to be the same, because the second transforms into the first under the 
change of integration variable $y=1/t$. However the separate integrals in 
(\ref{eq36}) are ill-defined because of a singularity at $t=1$. Nevertheless 
these integrals become well-defined and equal in the sense of Cauchy principal 
value. Therefore,
\begin{equation}
\phi(x)=P\int_0^\infty\frac{t^{-x}}{1-t}\,dt.
\label{eq37}
\end{equation}
To get rid of inconvenient principal value, we use quantum physicists' 
favorite formula (Sokhotski-Plemelj formula\footnote{Sokhotski-Plemelj formula
is  a relation between the generalized functions, that is it is assumed that
both sides of (\ref{eq38}) are multiplied by a smooth function, which is
non-singular in a neighborhood of the origin, then integrated over a range of 
$z$ containing the origin, and finally a limit $\epsilon\to 0$ is taken in
the results.}  
\cite{31,32})
\begin{equation}
\frac{1}{z\pm i\epsilon}=P\frac{1}{z}\mp i\pi\delta(z),
\label{eq38}
\end{equation}
which gives
\begin{equation}
\phi(x)\pm i\pi=\lim_{\epsilon\to 0}\int_0^\infty\frac{t^{-x}}
{1-t\mp i\epsilon}\,dt.
\label{eq39}
\end{equation}
Multiplying these two representations of $\phi(x)$, we get
\begin{equation}
\phi^2(x)+\pi^2=\lim_{\epsilon\to 0}\int_0^\infty\frac{t^{-x}}
{1-t+i\epsilon}\,dt\int_0^\infty\frac{s^{-x}}{1-s-i\epsilon}ds.
\label{eq41}
\end{equation}
We can substitute $s=y/t$ in the second integral, interchange the order
of integrations, solve the resulting simple integral in $t$,
\begin{equation}
\lim_{\epsilon\to 0}\int_0^\infty\frac{1}{(1-t+i\epsilon)(t-y-i\epsilon t)}
\,dt=-\frac{\ln{y}}{1-y},
\label{eq42}
\end{equation}
and end up with
\begin{equation}
\phi^2(x)+\pi^2=-\int_0^\infty\frac{y^{-x}\ln{y}}{1-y}\, dy.
\label{eq43}
\end{equation}
On the other hand, if we differentiate (\ref{eq37}) by $x$, we get  
\begin{equation}
\phi^\prime(x)\equiv\frac{d\phi(x)}{dx}=-\int_0^\infty\frac{t^{-x}\ln{t}}
{1-t}\, dt.
\label{eq44}
\end{equation}
(there is no longer a need for the principal value after differentiation, 
because the singularity softens and becomes integrable).
Comparing (\ref{eq43}) and (\ref{eq44}), we see that the function $\phi(x)$
satisfies differential equation
\begin{equation}
\phi^\prime(x)=\pi^2+\phi^2(x).
\label{eq45}
\end{equation}
Note that this differential equation is much simpler than the differential 
equation $F(x)F^{\prime\prime}(x)=(F^\prime(x))^2+F^4(x)$ obtained by Dedekind 
in \cite{28} for the function $F(x)=\Gamma(x)\Gamma(1-x)$.

Under the initial condition $\phi(1/2)=0$, which follows from the definition 
of $\phi(x)$, (\ref{eq45}) can be solved immediately:
\begin{equation}
x-\frac{1}{2}=\int_0^\phi\frac{d\phi}{\pi^2+\phi^2}=\frac{1}{\pi}\arctan{\frac
{\phi}{\pi}}.
\label{eq46}
\end{equation}
Therefore
\begin{equation}
\phi(x)=\pi\tan{\left(\pi x-\frac{\pi}{2}\right)}=-\pi\cot{\pi x},
\label{eq47}
\end{equation}
and this completes the proof of (\ref{eq33}).

\section*{Zeta function values at positive even integers}
The above approach can be easily generalized (this time without any
physics input) to enable a calculation of all $\zeta(2k)$.
Because of 
$$\sum\limits_{n=1}^\infty \frac{1}{n^{2k}}=\sum\limits_{n=1}^\infty \frac{1}
{(2n-1)^{2k}}+\frac{1}{2^{2k}}\sum\limits_{n=1}^\infty \frac{1}{n^{2k}},$$
we have
\begin{equation}
\zeta(2k)\equiv\sum\limits_{n=1}^\infty \frac{1}{n^{2k}}=\frac{2^{2k}}
{2^{2k}-1}\sum\limits_{n=1}^\infty \frac{1}{(2n-1)^{2k}}=\frac{1}{2^{2k}-1}
\sum\limits_{n=1}^\infty \frac{1}{(n-\frac{1}{2})^{2k}}.
\label{eq10}
\end{equation}
On the other hand, differentiating (\ref{eq7}) $2k-1$ times, we get
\begin{equation}
\psi_{2k-1}(1-x)=(2k-1)!\sum\limits_{n=1}^\infty \frac{1}{(n-x)^{2k}},
\label{eq11}
\end{equation}
where $\psi_n(x)=\frac{d^n}{dx^n}\psi(x)$ is the polygamma function. 
Therefore
\begin{equation}
\zeta(2k)=\frac{\psi_{2k-1}\left(\frac{1}{2}\right)}{(2^{2k}-1)(2k-1)!}.
\label{eq12}
\end{equation}
Since
$$\frac{\pi^2}{\sin^2{\pi x}}=-\pi\frac{d}{dx}\cot{\pi x},$$
differentiating (\ref{eq9}) $2k-)$ times, we get
\begin{equation}
\psi_{2k-1}(1-x)+\psi_{2k-1}(x)=-\pi\frac{d^{2k-1}}{dx^{2k-1}}\cot{\pi x}.
\label{eq13}
\end{equation}
It follows from this reflection formula that
\begin{equation}
\psi_{2k-1}\left(\frac{1}{2}\right)=-\frac{\pi}{2}\,s_{2k-1},
\label{eq14}
\end{equation}
where the $s_n$ numbers are defined through
\begin{equation}
s_n=\left . \frac{d^n}{dx^n}\cot{\pi x}\right |_{x=\frac{1}{2}}.
\label{eq15}
\end{equation}
Thanks to identity $\tan{\pi\left(x-\frac{1}{2}\right)}=-\cot{\pi x}$, we can 
express the $s_n$ numbers through more familiar tangent numbers \cite{11}
\begin{equation}
T_n=\left . \frac{d^n}{dx^n}\tan{x}\right |_{x=0}
\label{eq16}
\end{equation}
as $s_n=-\pi^n\,T_n$ and (\ref{eq12}) takes the form
\begin{equation}
\zeta(2k)=\frac{\pi^{2k}\,T_{2k-1}}{2(2^{2k}-1)(2k-1)!}.
\label{eq17}
\end{equation}
We can calculate tangent numbers recursively. To do so, note that
$$\frac{d^n}{dx^n}\tan{x}=\frac{d^{n-1}}{dx^{n-1}}\frac{1}
{\cos^2{x}}=\frac{d^{n-1}}{dx^{n-1}}\left(1+\tan^2{x}\right)=
\frac{d^{n-1}}{dx^{n-1}}\left(\tan{x}\,\tan{x}\right),$$
and apply the Leibniz formula for $n$-th derivative of a product of two 
functions. We get the recurrence relation
\begin{equation}
T_n=\sum\limits_{r=0}^{n-1}\binom{n-1}{r}T_r\,T_{n-1-r}.
\label{eq18}
\end{equation}
From the definition (\ref{eq16}) we find $T_0=0,\,T_1=1$ and it is not hard
to prove by induction that (\ref{eq18}) implies the vanishing of tangent 
numbers if their index is even. For an odd index, let's take $r=2m-1$ in  
(\ref{eq18}) to transform it into the form \cite{11}
\begin{equation}
T_{2k-1}=\sum\limits_{m=1}^{k-1}\binom{2k-2}{2m-1}\,T_{2m-1}\,T_{2(k-m)-1}.
\label{eq19}
\end{equation}
In principle (\ref{eq17}) and (\ref{eq19}) solve the generalized Basel 
problem and allow to calculate $\zeta(2k)$ at least for small values of $k$.
For example, we easily get $T_3=2,\,T_5=16,\,T_7=272,\,T_9=7936,\,
T_{11}=353792$ which imply
$$\zeta(4)=\frac{\pi^4}{90},\;\zeta(6)=\frac{\pi^6}{945},\;
\zeta(8)=\frac{\pi^8}{9450},\;\zeta(10)=\frac{\pi^{10}}{93555},\;
\zeta(12)=\frac{691\pi^{12}}{638512875}.$$
However for large values of $k$ more efficient algorithms are needed to
calculate tangent numbers. One of them can be found in \cite{12}.

A connection with Bernoulli numbers is established by the well-known formula
(see, for example, \cite{12}), valid for $n>1$,
$$B_n=-\frac{n\,T_{n-1}}{(2i)^n(2^n-1)}.$$
We also can use (\ref{eq12}), (\ref{eq14}) and
$$s_n=(-1)^n\frac{(2\pi)^{2n-1}(2^{2n}-1)}{n}\,B_n,$$
proved in \cite{13}. In either way we get the well known result
$$\zeta(2k)=(-1)^{k+1}\,\frac{(2\pi)^{2k}B_{2k}}{2\,(2k)!}.$$

\section*{Recurrence formula for $\zeta(2k)$}
With some extra effort, it is possible to obtain a nice recurrence 
relation for $\zeta(2k)$ \cite{14,15} which allows a calculation of $\zeta(2k)$
recursively, and thus also provides a solution of the generalized Basel 
problem.

Let's introduce another sequence of numbers related to the cotangent function:
\begin{equation}
S_n=\frac{d^n}{dx^n}\left .\left (x\,\cot{x}\right)\right|_{x=0}.
\label{eq20}
\end{equation}
In \cite{14} the coefficients of the Taylor expansion of $x\cot{x}$ were 
related to the values of $\zeta(2k)$. As an alternative to that argument, 
we will study the numbers $S_n$ and, using the insights from the previous 
section, relate them to $\zeta(2k)$.

Because of relation 
\begin{equation}
x\,\tan{x}=x\,\cot{x}-2x\,\cot{2x},
\label{eq21}
\end{equation}
which can be simply checked, the numbers introduced are related to the tangent 
numbers. Namely, differentiating (\ref{eq21}) $n$ times and setting $x=0$
produces the equality $n\,T_{n-1}=S_n-2^nS_n$. Therefore 
$$T_{2k-1}=-\frac{2^{2k}-1}{2k}\,S_{2k}$$ and
\begin{equation}
\zeta(2k)=-\frac{\pi^{2k}S_{2k}}{2(2k)!}.
\label{eq22}
\end{equation}
Now we obtain a recurrence relation for the $S_n$. By using 
$\cot^\prime{x}=-(1+\cot^2{x})$ it can be checked that \cite{14}
\begin{equation}
x\,(x\,\cot{x})^\prime=x\,\cot{x}-x^2\,\cot^2{x}-x^2.
\label{eq23}
\end{equation}
Differentiating both sides of this relation $n>2$ times and setting $x=0$,
we get
$$nS_n=S_n-\sum\limits_{r=0}^n\binom{n}{r}S_r\,S_{n-r}=
S_n-2S_n-\sum\limits_{r=1}^{n-1}\binom{n}{r}S_r\,S_{n-r}.$$
Therefore the $S_n$ numbers obey the recurrence relation 
\begin{equation}
(n+1)S_n=-\sum\limits_{r=1}^{n-1}\binom{n}{r}S_r\,S_{n-r}.
\label{eq24}
\end{equation}
As $S_0=1,\,S_1=0$, it follows from this recurrence relation (by induction)
that $S_n=0$, if $n$ is odd (if $n$ is odd, one of the numbers $r$, $n-r$ is 
also odd). Therefore, taking $r=2m$, we can write the recurrence relation 
(\ref{eq24}) in the form
\begin{equation}
(2k+1)S_{2k}=-\sum\limits_{m=1}^{k-1}\binom{2k}{2m}S_{2m}\,S_{2k-2m}.
\label{eq25}
\end{equation}
If we substitute (\ref{eq25}) into (\ref{eq22}), we get a recurrence relation
for the zeta-function that was called ``highly elegant'' in \cite{15}:
\begin{equation}
\left(k+\frac{1}{2}\right)\zeta(2k)=\sum\limits_{m=1}^{k-1}\zeta(2m)\,
\zeta(2k-2m).
\label{eq26}
\end{equation} 
In fact (\ref{eq26}) is equivalent to Euler's recurrence relation for 
Bernoulli numbers (independently found by Ramanujan \cite{16})
\begin{equation}
(2n+1)B_{2n}=-\sum\limits_{m=1}^{n-1}\binom{2n}{2m}B_{2m}\,B_{2n-2m}. 
\label{eq27}
\end{equation}  
Both (\ref{eq26}) and (\ref{eq27}) were rediscovered many times \cite{17}.
For example Williams \cite{18} thought he was the first to explicitly state 
the recurrence relation in the form (\ref{eq26}). Actually, this recurrence 
relation is given implicitly in Euler's work \cite{19} (paper E130 at
\url{http://www.math.dartmouth.edu/~euler/}), and is explicitly stated
at least as early as 1906 in the book \cite{20} (with the remark that
this recurrence relation is well known). Nowadays the proof of the recurrence 
relation (\ref{eq26}) is often given as an exercise in number theory courses
(see, for example, \cite{21}).

\section*{Acknowledgments}
The work is supported by the Ministry of Education and Science of the Russian 
Federation. The author thanks Professor Juan Arias de Reyna for indicating
several interesting references, as well as an anonymous referee for 
constructive comments which helped to improve the presentation.

Zurab K. Silagadze 

Budker Institute of Nuclear Physics and
Novosibirsk State University, 630 090, Novosibirsk, Russia.
silagadze@inp.nsk.su

\end{document}